\documentclass[11pt, oneside]{article}   	
\usepackage{geometry}                		
\usepackage{amssymb}
\usepackage{amsmath}
\usepackage{amsthm}
\usepackage{appendix}
\usepackage{graphicx}
\theoremstyle{definition}

\newtheorem{theorem}{Theorem}
\newtheorem*{theorem*}{Theorem}
\usepackage{bbm}
\newtheorem{lemma}[theorem]{Lemma}
\newtheorem*{remark}{Remark}
\newtheorem{corollary}{Corollary}

\title{ A Theoretical Analysis of the Stationarity of an Unrestricted Autoregression Process}
\author{Varsha S. Kulkarni 

\\
\\
Harvard University\\
}

\date{}							

\begin{document}
\maketitle
\begin{abstract}
The higher dimensional autoregressive models would describe some of the econometric processes relatively generically if they incorporate the heterogeneity in dependence on times. This paper analyzes the stationarity of an autoregressive process of dimension $k>1$ having a sequence of coefficients $\beta$ multiplied by successively increasing powers of $0<\delta<1$. The theorem gives the conditions of stationarity in simple relations between the coefficients and $k$ in terms of $\delta$. Computationally, the evidence of stationarity depends on the parameters. The choice of $\delta$ sets the bounds on $\beta$ and the number of time lags for prediction of the model.  $\footnote[1]{ Initial paper Acknowledgements: Thanks to Prabha Sharma, Ajit Iqbal Singh}$
\end{abstract}
\section{Introduction}

Autoregression models of time series represent stochastic dynamical systems. They are commonly used by econometricians as techniques to predict economic indicators based on the values in previous times. The randomness in the time evolution impacts the predictability and the number of lagged time steps included in the model. Some mathematical economists (Chiang and Wainwright, Im et al, Box et al, Judge et al, Samuelson)\cite{1,2,3,4,5} have investigated the time series properties of stationarity, ergodicity and causality theoretically to understand the dynamics of autoregressive models. In general, AR$(k)$ denotes an autoregressive model of dimension $k$. The research in this area has often addressed the theory of lower dimensional AR processes more than that of the higher dimensional models \cite{6}. In this paper, I analyze the stationarity conditions of a particular case of the unrestricted AR ($k$),  $k\geq 1$ model.

\medskip

\medskip
Previous research has focused on the stationarity of AR process with identical coefficients. In econometrics, there are, however, situations conforming to different kinds of AR$(k$) models. Consider the autoregressive model predicting the time evolution of a quantity or economic indicator when there are individuals or entities in the system showing varying levels of dependence on the previous values in time. For instance, the formation of inflationary expectations and the information propagation through the economy are heterogenous processes. 

 A prediction through an AR process may involve an unequal dependence on the different time lags, usually diminishing with increasing lags. The heterogeneity of the time dependence in this model maybe related to the asymmetry of information, inequality, or some other kind of phenomenon, and this may result in different models and their estimations. 

\medskip
This paper analyzes an unrestricted AR model wherein the coefficients of the first two time lags are equal, and the coefficients of the times beyond decrease sequentially with the number of previous time steps taken.

\section{Model}
The model is given as
\medskip

\begin{equation}
\pi(t) = \sum_{i=1}^{k} b_i \pi(t-i) +u(t)   ,k\geq 1
\end{equation}
where $u$ denotes the innovation term such that it has mean $Eu=0$ and variance as $\sigma^2$. And $\exists$ small $\epsilon_1, \epsilon_2>0$ so that $|b_2-b_1| = \epsilon_1$. Choosing a suitable $\epsilon_1$ and $0<\delta<1$ generalizes $b_i=\delta b_{i-1},   \forall i > 2$ . Also $\delta b_{i-1} >\epsilon_2$,  $\forall i>2$. 

\vspace{1.5\baselineskip}

The system in Eq.(1) describes a situation when there can be any number of past time steps or lags considered in the model with respect to the current time so that the $k^{th}$ term of the series would be at least equal to a chosen arbitrarily small number $\epsilon_2>0$. This restricts the permissible number of terms $k$ in the model, however, given that $\delta<1$ by choice, the variation in the number of terms can be very large. 

\vspace{1.5\baselineskip}

In a $k$ dimensional system, one may associate the coefficients with the proportions of individuals or firms in a population or economy according to their dependence or the relative importance placed on the values of the economic time series $\pi(t)$ in the previous times. In this context, I consider the first two $(b_1,b_2)$ coefficients denoting two different categories as being equal whereas all the other coefficients as proportionate to them.

The system in Eq(1) is equivalent to 
\begin{equation}
\pi(t)+\beta\pi(t-1)+\beta\pi(t-2)+\beta\delta\pi(t-3)+\beta\delta^2\pi(t-4)\cdots \beta\delta^{k-2}\pi(t-k) =u(t)
\end{equation}

where $b_i=-\beta$, $\forall i=1,2$ and $b_i=-\beta\delta^{i-2}$ $\forall i >2$ .

In this paper, I find the bound on the coefficient $\beta$ that satisfies the condition to achieve convergence to stationarity of the unrestricted AR$(k)$ system defined in Eq.(2). Researchers have solved this for AR$(k)$ process with equal coefficients (or weights), but the case presented in this paper is unique in the sense that the coefficients in the system bear a sequential relation.  Every coefficient is discounted in proportion to the coefficient at previous time step. 

\section{Analysis and proof}
The rest of this paper analyses the system and provides the proof. It depends on some well known results and theorems. 

The AR$(k)$ process in Eq.(1) has a characteristic polynomial of degree $k$ associated with it and (in accordance with algebraic theorem) every such polynomial has $k$ roots. 

The conditions on stationarity of AR ($k$) process defined above in (1) are equivalent to those of the  convergence of a polynomial (Box et al, Judge et al, Samuelson, Carmichael). From Eq.1
\begin{equation}
\pi (t) (1-b_1 x – b_2 x^2 -\cdots - b_k x^k) = b(x) \pi(t)
\end{equation}

In particular, the unrestricted AR process is stationary if the polynomial $b(x)=0$ polynomial has roots beyond the unit circle. This implies the condition that the roots of
\begin{equation}
 a^k – b_1 a^{k-1} -b_2 a^{k-2}….. -b_k=0
 \end{equation}
  would be lower than 1. In other words, the roots of (3) should be lower than 1 for AR($k$) process to converge to stationarity.
\vspace{1.5\baselineskip}

Researchers in economics and mathematics have shown that in order to ascertain the convergence of temporal paths, the conditions are derivable without the estimation of the roots of the characteristic equation. They have achieved this by utilizing some of the theorems including those of algebra, Rouche and Hurwitz theorems in complex analysis and Cauchy resulting in bounds and relations between roots and coefficients (Chiang and Wainwright, Samuelson, Carmichael)\cite{1,5,7}. In accordance with the Schur theorem, it is necessary and sufficient for the $k^{th}$ degree polynomial to have roots lower than unity, that all of its $k$ determinants are positive. 
Hence in this case, we rewrite the polynomial in Eq. (4) 
\begin{equation}
a^k+\beta a^{k-1} +\beta a^{k-2} +\beta\delta a^{k-3} +\beta\delta^{2} a^{k-4} +\cdots+\beta\delta^{k-2}
\end{equation}
This polynomial in Eq. (5) has coefficients as a geometric sequence. $A=A_m$ is a matrix of the structure of the blocks 
\begin{equation*}
A = 
\begin{pmatrix}
c & d^{\prime} \\
d & c^{\prime}
\end{pmatrix}
\end{equation*}
These blocks consist of matrices $c$ and $d$ and their transposes, which can be identified in A. In this context, the AR$(k)$ process is stationary if and only if $|A| >0$,  $\forall m=1,2..k$.

\medskip

$A$ is a square matrix of dimension $2m\times 2m$ matrix characterizing the polynomial wherein the coefficients $\beta$ are multiplied by the sequence of constants $\delta^0, \delta^0, \delta^1,\delta^2,\delta^3\cdots \delta^{k-2}$.

\begin{equation*}
A = 
\begin{pmatrix}
1 & 0 & 0& \cdots & 0  &   \beta\delta^{k-2} &\beta\delta^{k-3}  & \beta\delta^{k-4}& \cdots & \beta \\
\beta & 1 & \cdots &0 &  0 &   0 & \beta\delta^{k-2} &\beta\delta^{k-3} & \cdots & \beta \\
\vdots  & \vdots   & \vdots & \vdots  &\vdots & \vdots  & \vdots   & \vdots & \vdots &\vdots  \\
\vdots  & \vdots   & \vdots & \vdots  \ddots &\vdots & \vdots  & \vdots   & \vdots & \vdots  \ddots &\vdots \\
\beta\delta^{k-3} & \beta\delta^{k-4} & \cdots & \beta & 1 &  0 & 0 & \cdots & 0 & \beta\delta^{k-2}\\
\\
\beta\delta^{k-2}  & 0 & 0 &\cdots&0 &  1 & \beta & \beta &\cdots & \beta\delta^{k-3}\\
\beta\delta^{k-3} & \beta\delta^{k-2} & 0 & \cdots & 0 & 0 & 1 & \beta & \cdots & \beta\delta^{k-4}\\
\vdots  & \vdots   & \vdots & \vdots  &\vdots & \vdots  & \vdots   & \vdots & \vdots &\vdots  \\
\vdots  & \vdots   & \vdots & \vdots  \ddots &\vdots & \vdots  & \vdots   & \vdots & \vdots  \ddots &\vdots \\
\beta & \beta & \delta\beta & \cdots & \beta\delta^{k-2} & 0 & 0 & 0 & \cdots & 1
\end{pmatrix}
\end{equation*}

\vspace{1.5\baselineskip}

Evidently, if $\delta=1$ then the matrix $A$ represents the case of equal coefficients that has been studied previously (Im et al) \cite{2}. The spectral analysis of that system yields multiple eigenvalues of same magnitude indicating the inherent instability particularly for large $k$. Here, $0<\delta<1$, and the sequence of constants alters the matrix of coefficients and its properties while rendering the analysis of the matrix intricate. The triangularization of $A$ through Schur technique is algorithmic which may require a lot of numerical analysis, particularly as the size of the matrix becomes greater than 3. However, one may observe that as the dimension increases, the matrix A tends to get sparser. In addition, since $\delta$ is small, it creates a sparser matrix $A$. This has a simpler approximation if obtained through perturbation of the matrix and its properties. It gives a circle around eigenvalues as a range of stability and the deviation depends on the size and nature of the perturbation. 

\medskip
In accordance with the procedure for obtaining eigenvalues delineated earlier (Im et al)\cite{2}, $A$ can be written in terms of a square block matrix $v$, identity matrix $J$ and $\Delta$ matrix obtained by taking the Hadamard product. If $v$ is defined as $v = \left(\begin{smallmatrix} L_m & L_m^{\prime}\\ L_m & L_m^{\prime} \end{smallmatrix}\right)$, taking $L_m$ as a $m\times m$ lower triangular matrix such that $L(i,j) =1$ if 

\medskip
 $i\geq j$ and $L(i,j) =0$ otherwise.

\vspace{1.5\baselineskip}

Then 
$A=[v \beta + J(1-\beta)] \circ \Delta$. $\Delta_{ij}=A_{ij}$ whenever $A_{ij}=1$ and $A_{ij}=0$, and $\Delta_{ij}=f(\delta)$  otherwise. $\Delta$ is a square matrix of same dimension $2m$ such that the elements of $\Delta$ are identical to those of $A$ but without the coefficients $\beta$ in the elements. The matrix $v$ is a matrix of same dimension as A ($=2m$) and $v_1$ is a $2\times 2$ matrix with all elements as $1$.  $\circ$ is the Hadamard product. The eigenvalues of A, therefore, are bounded and dependent on the product of these matrices. (If all elements of $\Delta$ are 1, then it reduces to the all equal coefficients system studied previously). However, this is merely a slight perturbation to the system and therefore the eigenvalues should not be very different from those of $v \beta + J(1-\beta)$ and I incorporate this using an additive perturbed matrix, $R$ instead. $R$ is defined as

\begin{equation*}
R= 
\begin{pmatrix}
a_m^1 & a_m^{2\prime} \\
a_m^2 & a_m^{1\prime}
\end{pmatrix}
\end{equation*}

where
\medskip

\begin{equation*}
a_{m}^1 = 
\begin{pmatrix}
0 & 0 & 0  & \cdots & 0\\
\delta^0-1 & 0 & 0 & \cdots & 0\\
\delta^0-1 & \delta^0-1 & \delta^0-1 & 0 & \cdots & 0\\
\delta-1 & \delta^0-1 & \delta^0-1 & \cdots & 0\\
\delta^2-1 & \delta-1 &\delta^0-1 & \cdots & 0\\
\vdots  & \vdots  & \ddots & \vdots & 0  \\
\delta^{k-3}-1 & \delta^{k-4}-1 & \cdots &0 & 0 
\end{pmatrix}
\end{equation*}

and 

\begin{equation*}
a_{m}^2 = 
\begin{pmatrix}
\delta^{k-2}-1 & 0 & 0 & 0 & \cdots & 0\\
\delta^{k-3}-1 & \delta^{k-2}-1 & 0 & 0 & \cdots & 0\\
\vdots  & \vdots  & \ddots & \vdots  \\
\delta-1 & \delta-1 & \delta^2-1 & 0 & \cdots & 0\\
0 & \delta-1 & \delta-1 & 0 & \cdots & \delta^{k-2}-1\\
\end{pmatrix}
\end{equation*}
\medskip

 The original matrix can now be expressed as
 
 $A=[v +R]\beta + J(1-\beta)$\\

or
\begin{equation}
A=V\beta + J(1-\beta)
\end{equation}
taking $V=v +R$\\

Any vector can be an eigenvector of the Identity matrix with eigenvalue 1. Thus the eigenvalues of A are equal to $\lambda(V) (\beta) +(1-\beta)$ where $\lambda(V)$ denotes the eigenvalues of $V$. According to a well known Bauer-Fike theorem, the eigenvalues of $V$ will be inside the circle surrounded by eigenvalues of $v$ of radius as constant $\kappa$. In particular, \\

\medskip
$\lambda(V) \leq \lambda(v)+ \kappa_m $ and $\lambda(V) \geq \kappa_m- \lambda(v),  $where $\kappa_m=KB$ is the constant defined as

\medskip
 $B=B_m=||X||^{-1}||X||$ and $K=K_m=||R||$, $X$ being the matrix of eigenvectors of $v$. 
 
 \medskip
 \begin{remark}
 The number of terms of previous times in the AR$(k)$ model is bounded. 
 \medskip
 \end{remark}
 The remark statement is consistent with the definition and specification of parameters. In general, the $n^{th}$ term of the sequence converges as $\exists \epsilon>0$ such that $\beta\delta^n \approx \epsilon$. This is a Cauchy sequence and so the series can be truncated after some of the previous times. In other words, it specifies the starting point of the model. Therefore, $n$ is constrained by $\beta$, $\delta$ . It means that the number of terms and the dependence on the previous time steps diminishes depending on $\epsilon$ and other parameters.
 
 \vspace{1.5\baselineskip}

\begin{lemma}  $K_m$ is finite

\medskip
Proof:\\
$K_m$ is the norm of the matrix $R$ and we show here that it is bounded. The norm is\\
 $||R|| = \sqrt{\sum_{i,j}(R_{ij}^2)}$\\
This becomes upon simplification 
 \vspace{1.5\baselineskip}
$K=K_m = \\
\sqrt{2m[(\delta^{k-2}-1)^2+ (\delta^{k-3}-1)^2+(\delta^{k-4}-1)^2+\cdots+(\delta^3-1)^2+(\delta^2-1)^2+2(\delta-1)^2]}$.

 \vspace{1.5\baselineskip}
This expands to \\

\medskip
$K=\sqrt(2m[\delta^{2(k-2))}+\delta^{2(k-3))}+\cdots+\delta^6+\delta^4-2[\delta^{k-2}+\delta^{k-3}+\cdots+\delta^3+\delta^2]+m-3+2(\delta-1)^2])$
 \vspace{1.5\baselineskip}
and converges to\\

\medskip
$\sqrt((2m[\frac{\delta^4(1-\delta^{2(m-3))}}{1-\delta^2}]-2\frac{\delta^2(1-\delta^{m-3})}{1-\delta}+(m-3)+2(\delta-1)^2])$\\

 \medskip
Thus, 
\begin{equation}
K=\sqrt{(2m\left(\frac{\delta^2(1-\delta^{k-3}}{1-\delta}\{\frac{\delta^2(1+\delta^{k-3})}{1+\delta}-2\}+m-3+2(\delta-1)^2\right)}
\end{equation}

\vspace{1.5\baselineskip}

It must be noted that, $K_1 = \sqrt{2m(\delta^{k-2}-1)^2))}$, $K_2=\sqrt{2m(\delta^{k-2}-1)^2+(\delta^{k-3}-1)^2)}$ and so on. 

The polynomial in the square brackets under square root sign in Eq. (7) depends on constant functions $\delta$ and the parameter $m$ nonlinearly (increasing). The behavior of $m$ dominates that of the constant functions of $\delta$ and, therefore, $K$. It is evident that $min (\kappa_1, \kappa_{m>1}) = \kappa_1$. At fixed levels of $\beta$, $\delta$, the value of $K$ depends only on $m$. But $k$ is constrained to be around a certain value set by the choice of $0<\beta<1$ and $0<\delta <1$. Hence, $K_m$ remains finite within the realm of this context. Thus according to the lemma, this indicates that the eigenvalues are constrained to a finite range of stability.

\end{lemma}

\vspace{1.5\baselineskip}

\begin{lemma}  The eigenvalues of $A$, $\lambda_i (A)$ satisfy
\[
\begin{cases}
(m-\kappa_m)\beta +1  \leq  \lambda_i (A)\leq (m+\kappa_m)\beta+1  & i=1\\
1-\kappa_m\beta  \leq \lambda_i (A)  \leq  1+ \kappa_m\beta & i=2,3,..m\\
1-(\kappa_m+1)\beta  \leq \lambda_i (A) \leq  1+(\kappa_m-1)\beta & i=m+1,\cdots 2m
\end{cases}
\]
\end{lemma}

Proof\\

To find the eigenvalues of $v$, we can use a similarity transformation as shown earlier in the case of equal coefficients \cite{2}. According to this transformation, $C=\zeta^{-1} v \zeta$ with the matrix $\zeta$ defined as \begin{equation*}
\zeta= 
\begin{pmatrix}
I(k) & 0(k) \\
I(k) & I(k)
\end{pmatrix}
\end{equation*}
 and \begin{equation*}
\zeta^{-1}= 
\begin{pmatrix}
I(k) & 0(k) \\
-I(k) & I(k)
\end{pmatrix}
\end{equation*}

Therefore \begin{equation*}
C= 
\begin{pmatrix}
L_m+L_m\prime & L_m\prime \\
0_m & 0_m
\end{pmatrix}
\end{equation*}
 
 and eigenvalues of $v$ are same as those of $C$ which means, $\lambda(v)=\lambda(J)+\lambda(\mathbbm{1})$. Here, $\mathbbm{1}_{i,j}=1 \forall i,j$. \\
 
 \medskip
 Now $\lambda_i(\mathbbm{1}_m)=m$ for $i=1$ and $\lambda_i(\mathbbm{1}_m)=0$ for $i=2,3..m$. So, \\
 
 \medskip
 $\lambda_i(v)=m+1$ for $i=1$, $\lambda_i(v)=1$ for $i=2,3..m$, and $\lambda_i (v)=0$ for $i=m+1 \cdots 2m$
 
  First consider the case of lower than or equal to inequality which yields
\begin{equation*}
\lambda_i (A) \leq \begin{cases}
m\beta+1+ \kappa_m \beta & i=1\\
1+ \kappa_m\beta & i=2,3,..m\\
1-\beta + \kappa_m \beta & i=m+1,\cdots 2m.
\end{cases}
\end{equation*}

Then the second case of greater than or equal to as

\begin{equation*}
\lambda_i (A) \geq \begin{cases}
m\beta+1- \kappa_m \beta & i=1\\
1- \kappa_m\beta & i=2,3,..m\\
1-\beta - \kappa_m \beta & i=m+1,\cdots 2m.
\end{cases}
\end{equation*}

\vspace{1.5\baselineskip}

\begin{theorem*}
A necessary and sufficient condition that the time evolution of the autoregressive model \\
$\pi(t)+\beta\pi(t-1)+\beta\pi(t-2)+\beta\delta\pi(t-3)+\beta\delta^2\pi(t-4)\cdots \beta\delta^{k-2}\pi(t-k)$, \hskip 1 em $k>1$, $\beta>0$ and 

\medskip
$0<\delta<1$, converges to stationarity is that $0<\beta< \frac{1}{\kappa_k(\delta)+1}$.
\end{theorem*}
\medskip
Proof\\ 
The condition for stationarity is that the determinant $|A_m|>0$. Hence, the determinant as expressed in this equation is given as\\
\begin{equation}
|A_m|=\prod_{j=1}^{2m} \lambda_j (A_j)    \hskip 2em      m=1,2,..k
\end{equation}

The results of lemmas 1 and 2 specify the eigenvalues of $A$ and the range depending on the dimension considered. Therefore, 
\begin{equation}
((m-\kappa_m)\beta+1)(1-\kappa_m\beta)^{m-1}(1-(\kappa_m+1)\beta)^m \leq |A_m| \leq ((m+\kappa_m)\beta+1)(1+\kappa_m\beta)^{m-1}((\kappa_m-1)\beta+1)^m
\end{equation}
In the first case of lower than inequality

\medskip

$|A_m| \leq (m\beta+\kappa_m\beta+1)(1+\kappa_m\beta)^{m-1}(\kappa_m\beta+1-\beta)^m$.

\medskip
It requires that $ (m\beta+\kappa_m\beta+1)(1+\kappa_m\beta)^{m-1}(\kappa_m\beta+1-\beta)^m >0$. This results in \\

\begin{equation}
\beta>-\frac{1}{m+\kappa_m}.
\end{equation}
\vspace{1.5\baselineskip}

Then in the case of greater than inequality, $\lambda (V)>\lambda(v)$, the lower bound is
\vspace{1.5\baselineskip}
$|A_m| \geq (m\beta-\kappa_m\beta+1)(1-\kappa_m\beta)^{m-1}(1-\kappa_m\beta-\beta)^m$
The condition of stationarity now requires that $\beta$ be chosen so that the right side of the inequality is at least greater than 0. This means 
\begin{equation}
\beta<\frac{1}{\kappa_m+1}
\end{equation}

\medskip

Thus from Eqs. (10) and (11), the condition is $\beta<\frac{1}{\kappa_m+1}$ and $\beta>-\frac{1}{m+\kappa_m}$. $\kappa_m$ increases with $m$, and so $-\frac{1}{k+\kappa_k(\delta)}<\beta< \frac{1}{\kappa_k(\delta)+1}$. If $\beta>0$, it simply reduces to $0<\beta< \frac{1}{\kappa_k(\delta)+1}$.

\medskip
 Thus it is interesting that the stationarity of this unrestricted AR$(k)$ process is given by a simple relation, determined by the coefficients and $k$. \\
 The attempt to reproduce this result on a time series in the appendix shows that the bounds on $\beta$ for obtaining stationarity may not be as strict as given in the theorem. Rather, it manifests clearly in terms of the variation of the strengths of evidence of stationarity within versus outside these bounds. In other words, to obtain a perceptible inference of non stationarity, one should take $\beta$ much greater than specified by this boundary.

\begin{corollary} 
The stationary of this $AR(k)$, $\hskip 0.025 em k>1$ model requires sum of all its coefficients to be lower than 1.

\medskip
Proof:\\
The theorem specifies the allowed interval around $\beta$ when the AR$(k)$ model is . Therefore, here, the sum of coefficients is \\
$\beta+\beta(\delta^0+\delta^1+\delta^2+\cdots \delta^{k-2})$. The sum converges and 
\begin{equation}
\beta+\beta \left [ \frac{1-\delta^{k-1}}{1-\delta}\right ] < \frac{1}{\kappa_k+1}\left (\frac{2-\delta-\delta^{k-1}}{1-\delta}\right)
\end{equation}

Then for the sum to be $< 1$, I show that ${\kappa_k+1}>\left (\frac{2-\delta-\delta^{k-1}}{1-\delta}\right)$. By the construction, $\kappa_k=K_k B_k >>\kappa_1$ when $k>>1$. 

\medskip
And, $\kappa_1 = K_1$ because $v_1$ has a unitary matrix of eigenvectors as $x= \left(\begin{smallmatrix} -\sqrt(2)&\sqrt(2)\\ \sqrt(2)&\sqrt(2) \end{smallmatrix}\right)$, therefore, $B_1=1$ and $\kappa_1=\sqrt(2) (1-\delta^{k-2})$. If we consider the case of a relatively small $\delta$, then ${\kappa_1+1}>\left (\frac{2-\delta-\delta^{k-1}}{1-\delta}\right)$. In particular if $\delta\approx 0$, for the stationary AR(2) process, $\beta<\frac{1}{\kappa_2(\delta)+1}<1$. In case of larger value of $\delta$, both $K_k$ and $B_k$ increase with $k$ or as $k\to\infty$. $B_k$ increases as a function of $k$ and $k^2$, so ${\kappa_k+1}>\left (\frac{2-\delta-\delta^{k-1}}{1-\delta}\right)$
and this means that the sum total of the coefficients of the unrestricted AR model is $< 1$. 
\end{corollary}

\begin{corollary}
In the stationary AR$(k)$ model, $k>1$ the lower bound on $\delta$ is dependent on $k$ as $\delta> \{1+\frac{1}{\sqrt{2}} \left(\frac{\beta-1}{\beta}\right)\}^{\frac{1}{k-2}}$.
\end{corollary}

Proof:\\
We have the AR($k$) model, $k>1$. If this process converges to stationarity, then for a fixed $\beta$ as constant, according to the theorem, $\beta< \frac{1}{\kappa_k(\delta)+1}$ and that $\kappa_1+1<\frac{1}{\beta}$. Taking $\kappa_1=\sqrt(2) (1-\delta^{k-2})$ from the Corollary 1, we have 

\medskip
$\delta> \{1+\frac{1}{\sqrt{2}} \left(\frac{\beta-1}{\beta}\right)\}^{\frac{1}{k-2}}$. This establishes the lower bound on $\delta$. It is consistent with $\delta<1$. It is interesting as it indicates that it decreases as $k$ increases. Thus the choice of $\delta$ restricts the number of lagged terms in the model. The example of the time series and simulation in the Appendix also demonstrates this.

\ \section*{Appendix} 

Figure 1 below shows the time series of Bombay Stock Exchange (BSE) market index, called the BSE Sensex. The coefficient of the AR(1) model of the first differenced time series $\pi(t)$ for the first 100 data points (days) estimated through linear regression is less than 1, which implies stationarity of the model. Then the coefficients of the generic model in Eqs. (1,2) computed using linear regression for $k=5$ are used to fit the model for another time period. The Akaike Information Critierion is usually applied to know the number of lags though here we take the coefficients sequentially so that greater time lags terms would diminish gradually. The estimated regression prediction on the next 40 days is not significant, however, it shows reasonably accuracy. Further, the simulated time series using the regression coefficients is non stationary and does not predict significantly.  

\begin{figure}[ht]
  \begin{center}
  \includegraphics[width=16cm, height=7cm]{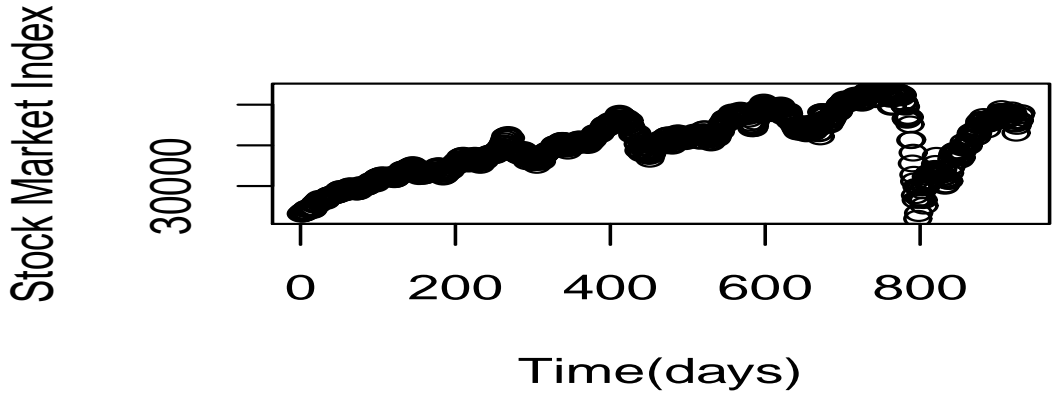}
  \end{center}
  \caption{Time evolution of the BSE Sensex, $\Delta t =1$ day, from January 2017 to December 2020}
  \label{states}
\end{figure}

\begin{figure}[ht]
  \begin{center}
  \includegraphics[width=16cm, height=7cm]{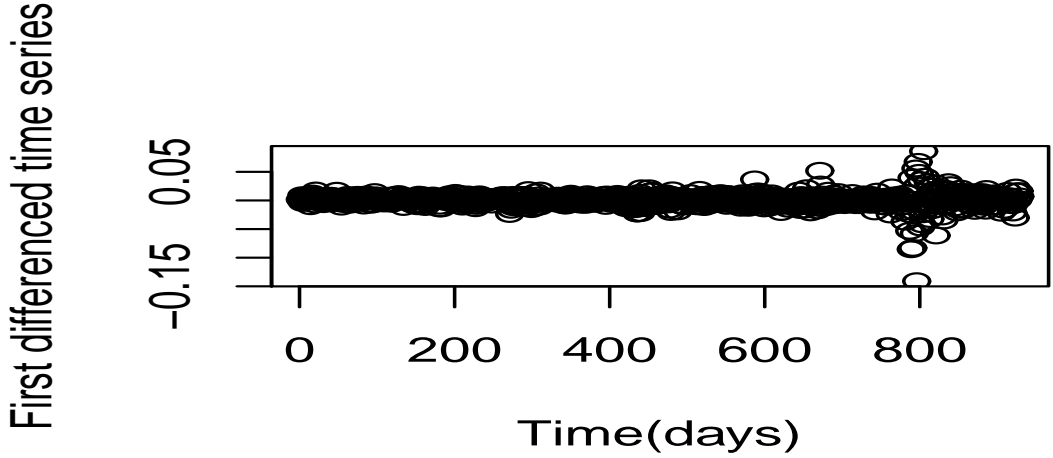}
  \end{center}
  \caption{The time series of first differences $\pi(t+1)-\pi(t)$}
\end{figure}

\medskip

In accordance with the theorem, the choice of $\delta$ determines $\beta$ for attaining stationarity. The model in Eqs.(1,2) simulated for $k=5$ shows varying strengths of evidence for stationarity (Table 1) based on the Augmented Dickey Fuller test statistic. \\
$\pi(t)+\beta\pi(t-1)+\beta\pi(t-2)+\beta\delta\pi(t-3)+\beta\delta^2\pi(t-4)+\beta\delta^3\pi(t-5)$\\

\begin{table}[h]
\caption{Computed model stationarity}
\centering
\begin{tabular}{c c c c}
\hline\hline
$\beta$ & $\delta$ & Stationarity \\ [0.5ex] 
\hline
0.0001&0.5&Highly strong stationarity \\
0.3&0.5&Moderately stationary \\
0.35&0.5&Weak stationarity \\
0.4&0.5&Non stationarity \\
0.0001&0.55&Highly strong stationarity \\
0.3&0.55&Moderately stationary \\
0.35&0.55&Weak stationarity \\
0.4&0.55&Non stationarity \\
0.0001&0.6& Stationary \\
0.3&0.6&Moderately stationary \\
0.35&0.6&Non stationarity \\
0.4&0.6&Non stationarity \\
0.0001&0.65& Stationary \\
0.3&0.65&Moderately stationary \\
0.35&0.65&Non stationarity \\
0.4&0.65&Non stationarity \\
0.0001&0.7& Weak Stationary \\
0.3&0.7& Non stationary \\
0.35&0.7&Non stationarity \\
0.4&0.7&Non stationarity \\ [1ex]
\hline
\end{tabular}
\label{table:nonlin}
\end{table}

In the simulation, the initial values maybe taken as the consecutive data points from any time period within the time series for computing the model for about 50 time steps in future using different values of $\beta$. Table 1 shows that as $\beta$ increases beyond $ \frac{1}{\kappa_k(\delta)+1}$ ( Theorem), the evidence of stationarity weakens. However, it is very strong when $\beta$ is close to the interval set by the theorem. The simulated model also predicts the BSE Sensex with reasonable accuracy in the time period. The root mean square error for comparing the original data and prediction model: $\pi(t)+\beta\pi(t-1)+\beta\pi(t-2)+\beta\delta\pi(t-3)+\beta\delta^2\pi(t-4)+\beta\delta^3\pi(t-5)$ is slightly lower and hence the model computed for prediction is significantly more accurate than the one obtained by using the regression coefficients, an advantage of the model. The boundaries set by the theorem on $\beta$ are, however, relatively flexible on this time series.

\end{document}